\documentclass{amsart}
\newtheorem*{theorem}{Theorem}
\newtheorem*{lemma}{Lemma}
\newtheorem*{prop}{Proposition}

\newtheorem*{definition}{Definition}

\DeclareMathOperator\di{\nabla \cdot}
\DeclareMathOperator\cu{\nabla \times}

\DeclareMathOperator\conj{Conj}

\usepackage[retainorgcmds]{IEEEtrantools}
\begin{document}
\def\HCD{H(\mathrm{curl};\O)\boldsymbol{\cap} H(\mathrm{div};\O)}
\def\O{\Omega}
\def\p{\partial}
\def\R{\mathbb{R}}
\def\cP{\mathcal{P}}
\def\cE{\mathcal{E}}
\def\cH{\mathcal{H}}
\def\cN{\mathcal{N}}
\def\cM{\mathcal{M}}
\def\cT{\mathcal{T}}
\def\mZ{\mathbb{Z}}
\def\boe{{\bf e}}
\def\bu{{\boldsymbol u}}
\def\bv{{\boldsymbol v}}
\def\bn{\boldsymbol{n}}
\def\bzeta{\boldsymbol{\zeta}}
\def\curl{\nabla\times}
\def\div{\nabla\hspace{2pt}\cdot\hspace{2pt}}

\newcommand{\be}{\begin{equation}}
\newcommand{\ee}{\end{equation}}
\newcommand\sep{\; ; \;}

\def\<{\langle}
\def\>{\rangle}
\newcommand{\iref}[1]{(\ref{#1})}
\def\n{{\rm \bf n}}
\def\t{{\rm \bf t}}
\def\sq{\hfill $\diamond$\\}
\def\C{{\rm \hbox{C\kern -.5em {\raise .32ex \hbox{$\scriptscriptstyle
|$}}\kern-.22em{\raise .6ex \hbox{$\scriptscriptstyle |$}}\kern .4em}}}
\newcommand\trans{\mathrm T}

\def\Hilbert{P_{k-1}(\R^2)}
\def\ve{\varepsilon}
\def\perm{{\mathbb S}}
\def\nl{\newline}

\title[Nonconforming Element for
  ${H}(\text{curl};{\Omega})\cap{H}(\text{div};{\Omega})$]
  {Nonconforming Vector Finite\\ Elements for
  $\boldsymbol{H}(\text{curl};\boldsymbol{\Omega})\cap
  \boldsymbol{H}(\text{div};\boldsymbol{\Omega})$}
\author{Jean-Marie Mirebeau}
%
%
%
\begin{abstract} 
We present a family of nonconforming vector finite elements of arbitrary order for
problems posed on the space $\HCD$, where $\O\subset\R^2$. 
This result was first stated as a conjecture by Brenner and Sung in \cite{BrennerConj}. In contrast an extension of the same conjecture to domains of $\R^3$ is disproved.
\end{abstract}
\maketitle
%
Let $\O$ be a domain of $\R^d$ where $d\in \{2,3\}$. 
As explained in \cite{BrennerConj} several problems involving the space $\HCD$, such as
 the cavity resonance problem and the
 acoustic fluid-structure interaction problem, can be solved using nonconforming finite element methods. In contrast conforming finite element methods cannot capture the solution of these problems under certain conditions. 

The accuracy of the approximate numerical solution of these problems can be improved if one uses finite elements which are not piecewise linear, but piecewise quadratic or of higher degree. For that purpose a quadratic nonconforming vector finite element for $\HCD$ was introduced in \cite{BrennerConj}, in the case of a two dimensional domain $\O\subset \R^2$. The paper \cite{BrennerConj} also contains a conjecture which suggests a way of constructing nonconforming vector finite elements of arbitrary degree $k$ for $\HCD$, for domains of $\R^2$ and of $\R^3$.

In order to state this conjecture and to formulate our results, we need to introduce some notations. We use boldfaced
 letters to represent vectors.  The space of polynomials of
 total degree $\leq k$ in $d$ variables is denoted by
 $P_k(\R^d)$, and the space of homogeneous harmonic polynomials
 of degree $k$ in $d$ variables is denoted by $\cH_k(\R^d)$.
For each $k\geq 1$ and $d\in \{2,3\}$ we define a space $\cP_{k,d}$ of vector fields on $\R^d$ as follows
\be
\label{defPkd}
\cP_{k,d} := [P_k(\R^d)]^d \oplus \left(\nabla \cH_{k+2}(\R^d) \oplus \cdots \oplus \nabla \cH_{2k}(\R^d)\right)
\ee
For any triangle $T$ if $d=2$ (resp. tetrahedron $T$ if $d=3$) we consider a collection $\cN_{k,d} = \cN_{k,d}^1 \cup \cN_{k,d}^2$ of linear functionals on $\cP_{k,d}$, which is defined as follows.
The elements of $\cN_{k,d}^1(T)$  define the moments on $T$ up to order $k-2$, for the $d$ components of the vector fields. Denoting by $\lambda_0, \cdots, \lambda_d$ the barycentric coordinates on the simplex $T$, and by $e_1, \cdots, e_d$ a fixed basis of $\R^d$, these functionals have the expression
$$
M_{i,\alpha}(v) := \int_T \lambda_0^{\alpha_0} \cdots \lambda_d^{\alpha_d} (\bv\cdot e_i),
$$
for all $1\leq i \leq d$ and all
$\alpha = (\alpha_0, \cdots, \alpha_d) \in \mZ_+^{d+1}$ such $\alpha_0+\cdots +\alpha_d=k-2$. The elements of $\cN_{k,d}^2(T)$  define the moments up to order $k-1$ on the $d+1$ edges (resp. faces) of $T$, defined by $F_j := \{z\in T :\ \lambda_j(z)=0\}$, $0\leq j\leq d$, again for the $d$ components of the vector fields. These functionals have the expression
$$
M_{i,j,\beta}(v) := \int_{F_j} \lambda_0^{\beta_0} \cdots \lambda_d^{\beta_d} (\bv\cdot e_i),
$$
for all $1\leq i \leq d$, all $0\leq j\leq d$, and all $\beta = (\beta_0, \cdots, \beta_d)\in \mZ_+^{d+1}$ such that $\beta_0+\cdots+\beta_d = k-1$ and $\beta_j = 0$. We used in this definition the convention $\lambda_j^{\beta_j} = 0^0 = 1$ on $F_j$.


Brenner and Sung formulated in \cite{BrennerConj} a series of conjectures, which depend on the two parameters $d\in \{2,3\}$ and $k\geq 1$.
\be
\label{MainConj}
\begin{array}{c}
\conj(k,d) : \text{ For any } T \text{ the elements of } \cP_{k,d} \text{ are uniquely determined}\\ 
\text{by the linear functionals in }\cN_{k,d}(T).
\end{array}
\ee
The conjectures $\conj(1,2)$ and $\conj(1,3)$ are true and correspond to the nonconforming Crouzeix-Raviart $P_1$ vector finite element. It was established in \cite{BrennerConj} that $\conj(2,2)$ is true, thus defining piecewise quadratic nonconforming vector finite elements in two space dimensions.

The purpose of this paper is to establish the following result : 
\begin{theorem}
For any $k\geq 3$ the conjecture $\conj(k,2)$ holds. In contrast the conjecture
$\conj(2,3)$ does not hold.
\end{theorem}
Our result therefore validates the construction of bi-dimensional vector finite elements of arbitrary degree proposed in \cite{BrennerConj}. On the contrary the three-dimensional quadratic vector finite element is invalid. Our result does not completely close the conjecture as the cases of three-dimensional vector finite elements of cubic or higher degree remain unsolved.


It was established in \cite{BrennerConj} that for all $d\in \{2,3\}$, all $k\geq 1$ and all $T$, one has 
$$
\dim \cP_{k,d} = \# \cN_{k,d}(T).
$$
Hence the conjecture $\conj(k,d)$ is equivalent to the following property : 
\be
\label{kernel}
\begin{array}{c}
\text{For all } T \text{ and all } \bv\in \cP_{k,d}, \\
 \text{if } \left( l(\bv) = 0 \text{ for all } l\in \cN_{k,d}(T) \right) \text{ then } \bv=0.
 \end{array}
\ee
In the first section of this paper we establish this property in the bi-dimensional case $d=2$ and for an arbitrary $k\geq 1$. In contrast we give in the second section a counter example in the three-dimensional case $d=3$ and $k=2$. 

\section{Proof of the bi-dimensional result}
 In this section the integer $k\geq 1$ is arbitrary but fixed. If $\bv = (v_1,v_2)\in \cP_{k,2}$ we remark that 
$$
 \curl\bv=\frac{\p v_2}{\p x_1}-\frac{\p v_1}{\p x_2}\in P_{k-1}(\R^2)\quad
 \text{and}\quad \div\bv=\frac{\p v_1}{\p x_1}+\frac{\p v_2}{\p x_2}\in P_{k-1}(\R^2).
$$
Our first lemma extends to degree $k$ an argument used in the initial paper \cite{BrennerConj}.
\begin{lemma}
Let $\bv\in \cP_{k,2}$. Let $T$ be a triangle and let us assume that 
$l(\bv) = 0$ for all $l\in \cN_{k,2}(T)$.
Then
\be
\label{vanishdicu}
\di \bv = \cu \bv = 0.
\ee
\end{lemma}

\proof
We first notice that $\di \bv$ and $\cu \bv$ are polynomials of degree $k-1$, and that the components of $\curl(\curl\bv)$ and of $\nabla(\div\bv)$ are polynomials of degree $k-2$.
In view of 
  Green's theorem and the vanishing moments of $\bv$, we have
 $$\int_T(\curl\bv)(\curl\bv)\,dx=\int_{\p T}(\bn\times\bv)(\curl\bv)\,ds
      +\int_T\bv\cdot\curl(\curl\bv)\,dx=0$$
 where $\bn$ is the outer unit normal along $\p T$.  Similarly, we have
 $$\int_T(\div\bv)(\div\bv)\,dx=\int_{\p T}(\bn\cdot\bv)(\div\bv)\,ds-
   \int_T\bv\cdot\nabla(\div\bv)\,dx=0.$$
The results follow.
\sq

We now rephrase the conjecture \iref{MainConj} in terms of complex functions, and for that purpose we introduce some definitions.
\begin{definition}
For any pair $\bv = (v_1,v_2)$ of real valued functions we define a complex valued function $P_\bv$ as follows
$$
P_\bv(x+iy) := v_1(x,y) -i  \, v_2(x,y) \text{ for all } (x,y)\in \R^2. 
$$
\end{definition}
We now note that the equations $\div \bv = \cu \bv =0$ are equivalent to the Cauchy-Riemann equations of $P_\bv$, 
namely 
$$
\frac{\partial \,\Re (P_\bv)}{\partial x} = \frac{\partial \, \Im (P_\bv)}{\partial y} \text{ and } \frac{\partial \, \Re (P_\bv)}{\partial y} = -\frac{\partial \, \Im (P_\bv)}{\partial x}
$$
where $\Re : \C\to \R$ and $\Im : \C \to \R$ respectively refer to the real and imaginary part. These equations characterize holomorphic functions.
Let us introduce for all $m\geq 1$  the space $\C_m$ of polynomials in the complex variable $z
=x+iy$ and of degree less or equal to $m$
$$
\C_m := \left\{P = \sum_{r=0}^m a_r z^r\sep (a_0, \cdots, a_r)\in \C^m\right\}.
$$
If $\bv \in \cP_{k,2}$ satisfies $l(\bv) = 0$ for all $l\in \cN_{k,2}$, then $P_\bv$ satisfies the Cauchy-Riemann equations  according to \iref{vanishdicu}, and therefore $P_\bv \in \C_{2k-1}$.

\begin{definition}
For any continuous function $P : \C\to \C$ and any $z_1, z_2\in \C$ we define 
\be
\label{defIP}
I_{z_1,z_2} (P) = \int_{t=0}^1 P(z_1 + t(z_2-z_1)) (z_2-z_1) dt = \int_S P(z) dz 
\ee
where $S\subset \C$ is the oriented segment from $z_1$ to $z_2$.
\end{definition}
Let $S$ be an edge of a triangle $T$ with endpoints $(x_1,y_1)$ and $(x_2,y_2)$ and let $z_1 = x_1+i y_1$ and $z_2 = x_2+i y_2$ be their complex coordinates. 
Let $\bv = (v_1,v_2)\in \cP_{k,2}$ be such that  $l(\bv) = 0$ for all $l\in \cN_{k,2}(T)$, and let $Q(x+ i y) := R_1(x,y) +i R_2(x,y)$ where $R_1,R_2 \in P_{k-1}(\R^2)$ are arbitrary. Since $\bv$ has vanishing moments up to order $k-1$ on the edges of $T$ we have
\be
\label{IPv}
I_{z_1,z_2}(P_\bv Q) = (z_2-z_1)\int_{t=0}^1 \big( (v_1 R_1 + v_2 R_2) + i (v_1 R_2 - v_2 R_1)\big)|_{(x(t),y(t))} dt  =0,
\ee
Where we used the notations $x(t) :=  x_1+t(x_2-x_1)$ and $y(t) :=  y_1+t(y_2-y_1)$.

We now define a bilinear form which is related to our conjecture.
\begin{definition}
For all $Z = (z_1,z_2,z_3)\in \C^3$ 
we define a bilinear form
$q_Z : \C_{2k-1} \times (\C_{k-1}\times \C_{k-1})\to \C$ as follows
$$
q_Z(P, (Q_1,Q_2)) := I_{z_1,z_2} (PQ_1) + I_{z_1,z_3} (PQ_2).
$$
\end{definition}
Let $T$ be a triangle and let $z_1=x_1+i y_1$, $z_2=x_2+iy_2$ and $z_3 = x_3+i y_3$ be the complex coordinates of the vertices of $T$.
If $\bv \in \cP_{k,2}$ is such that $l(\bv) = 0$ for all $l \in \cN_{k,2}(T)$ then $P_\bv \in \C_{2k-1}$ as previously noted. Furthermore, specializing \iref{IPv} to polynomials $Q\in \C_{k-1}$ we obtain
\be
\label{qPv}
q_Z(P_\bv, (Q_1,Q_2)) = 0 \text{ for all } (Q_1,Q_2)\in \C_{k-1}\times \C_{k-1}.
\ee
The purpose of the rest of this section is to show that the bilinear form $q_Z$ is nondegenerate. It then follows from \iref{qPv} that $P_\bv = 0$ and therefore that $\bv = 0$ which concludes the proof of the conjecture $\conj(k,2)$. 

We denote by $B :=  (1,z,\cdots, z^{2k-1})$ the canonical basis of $\C_{2k-1}$, and by 
$B^* := ( (1,0), (z,0) \cdots ,(z^{k-1},0), \ (0,1), \cdots , (0,z^{k-1}) )$ the canonical basis of $\C_{k-1} \times C_{k-1}$.
We denote by $M(Z)$, or $M(z_1,z_2,z_3)$, the matrix of $q_Z$ in the basis $B$ and $B^*$. Hence for all $1\leq i \leq 2k$ and all $1\leq j \leq k$ we have
$$
M(Z)_{i,j} = I_{z_1,z_2} (z^{i-1} z^{j-1}) \ \text{ and } \ M(Z)_{i,j+k} = I_{z_1,z_3} (z^{i-1} z^{j-1})
$$
It follows that 
\be
\label{defMz}
M(Z)_{i,j} = \frac{z_2^{i+j-1}-z_1^{i+j-1}} {i+j-1} \text{ and } M(Z)_{i,j+k} = \frac{z_3^{i+j-1}-z_1^{i+j-1}} {i+j-1} .
\ee
The explicit expression of the matrix $M(Z)=M(z_1,z_2,z_3)$, in the special case $k=2$, is given right after the end of the proof.

Our next proposition gives an explicit expression of $\det M(Z)$, therefore showing that $q_Z$ is non-degenerate. In the following $Z$ always refers to the triplet of complex variables $Z=(z_1,z_2,z_3)$.
\begin{prop}
One has
$$
\det M(Z) =  \alpha (z_1-z_2)^{k^2} (z_2-z_3)^{k^2} (z_3-z_1)^{k^2} 
$$
where $\alpha =  \frac{\left(\prod_{0\leq i \leq k-1} i!\right)^5}{\prod_{0\leq i\leq k-1} (2k+i)!}  > 0$. Therefore $q_Z$ is non-degenrate whenever $z_1$, $z_2$ and $z_3$ are pairwise distinct.
\end{prop}

\proof
We denote by $\perm$ the collection of all permutations $\sigma$ of the set $\{1,\cdots , 2k\}$, and by $\ve(\sigma)$ be the algebraic signature of such a permutation. We recall that 
\be
\label{detsigma}
\det M(Z) := \sum_{\sigma\in \perm} \ve(\sigma) \prod_{j=1}^{2k} M(Z)_{\sigma(j),j}. 
\ee
For any permutation $\sigma\in \perm$ one has 
$$
\sum_{j=1}^k ( j + \sigma(j) -1) + \sum_{j=1}^k ( j + \sigma(k+j) -1) = 3 k^2.
$$
It follows from \iref{defMz} that $\det M(Z) $ is a homogeneous polynomial in the variables $z_1,z_2,z_3$ and of degree $3 k^2$.
We also note for future use that 
\be
\label{sigmak2}
\sum_{j=1}^k ( j + \sigma(j) -1) \geq k^2
\ee
with equality if and only if $\sigma$ leaves invariant the sets $\{1,\cdots, k\}$ and $\{k+1, \cdots, 2k\}$.
For any $c\in \C$ we define two $2k \times 2k$ triangular matrices $P(c)$ and $P^*(c)$ associated with the following changes of basis on $\C_{2k-1}$ and $\C_{k-1}\times \C_{k-1}$ respectively
\begin{eqnarray*}
P(c) B &=& (\ 1,\ z+c, \cdots, \ (z+c)^{2k-1})\\
P^*(c) B^* &=& (\ (1,0),\cdots,\  ((z+c)^{k-1},0),\  (0,1), \cdots,\ (0, (z+c)^{k-1})\,)
\end{eqnarray*}
One easily sees that the matrices $P(c)$ and $P^*(c)$ are upper-triangular and have ones on the diagonal, hence $\det P(c) = \det P^*(c) = 1$. 

Since 
$$
I_{z_1+c,z_2+c}(z^iz^j) = I_{z_1,z_2}((z+c)^i (z+c)^j)
$$ 
we obtain 
$$
M(z_1+c,\ z_2+c,\ z_3+c) = P(c)^\trans M(Z) P^*(c). 
$$
Recalling that $\det P(c) = \det P^*(c) = 1$, and choosing $c=-z_1$, we obtain 
$$
\det M(0,\ z_2-z_1,\ z_3-z_1) = \det M(Z). 
$$
The explicit expression of the matrix $M(0,\ z_2-z_1,\ z_3-z_1)$, in the special case $k=2$, is given right after the end of the proof.
It follows from \iref{defMz}, \iref{detsigma} and \iref{sigmak2} that 
the polynomial $\det M(Z)$ is a multiple of $(z_2-z_1)^{k^2}$. Similarly, $\det M(Z)$ is a multiple of $(z_3-z_1)^{k^2}$.

Substracting column $k+i$ from column $i$, for all $1\leq i\leq k$, we find that 
$$
\det M(z_3, z_2, z_1) = (-1)^k \det M(Z) 
$$
and therefore 
$\det M(Z)$ is also a multiple of $(z_3-z_2)^{k^2}$. Since $\det M(Z)$ is a polynomial of degree $3k^2$ in the complex variables $z_1,z_2, z_3$, and since $(z_1-z_2)^{k^2}$, $(z_2-z_3)^{k^2}$ and $(z_3-z_1)^{k^2}$ have no common factors, there exists a constant $\alpha \in \C$ such that 
$$
\det M(Z) = \alpha \ (z_1-z_2)^{k^2} (z_2-z_3)^{k^2} (z_3-z_1)^{k^2}.
$$
In order to compute the constant $\alpha$, and to show that $\alpha\neq 0$, we remark that it is the coefficient of  $z^{k^2}$ in the polynomial $\det M(0,z,1)= \alpha (-z)^{k^2} (z-1)^{k^2}$. The explicit expression of this matrix, in the special case $k=2$, is given right after the end of the proof.

The contribution of a permutation $\sigma\in \perm$ to $\det M(0,z,1)$ is a monomial which has degree $k^2$ if and only if \iref{sigmak2} is an equality. 
Denoting by $\perm^*$ the collection of permutations of the set $\{1, \cdots, k\}$, we obtain that 
$\det M(0,z,1)$ equals
$$
\left(\sum_{\sigma_1 \in \perm^*} \ve(\sigma_1) \prod_{j=1}^k M(0,z,1)_{j,\sigma_1(j)} \right) \left(\sum_{\sigma_2 \in \perm^*} \ve(\sigma_2) \prod_{j=1}^k M(0,z,1)_{j+k,\, \sigma_2(j)+k} \right)+ \mathcal O(z^{k^2+1}).
$$
Hence using \iref{defMz}
$$
\det M(0,z,1) = z^{k^2} \det\left(\frac 1 {i+j-1}\right)_{1\leq i,j\leq k} \det \left(\frac 1 {i+j+k -1}\right)_{1\leq i,j\leq k} + \mathcal O(z^{k^2+1})
$$
This expression gives the value of $\alpha$ as the product of two Cauchy determinants, which can be computed using the formula, established in \cite{CauchyDet} \S I.1.3,
$$
\det \left(\frac 1 {a_i+ b_j}\right)_{1\leq i,j\leq k} = \frac{\underset{1\leq i<j \leq k}\prod (a_i - a_j) \underset{1\leq i<j \leq k}\prod (b_i - b_j)}{\underset{1\leq i,j\leq k} \prod (a_i+b_j)}. 
$$
This concludes the computation of $\det M(Z)$.
\sq

The following explicit expressions of the matrix $M$, if $k=2$ and $z,z_1,z_2,z_3\in \C$, may be useful to the reader
\begin{eqnarray*}
\label{exMz}
M(z_1,z_2,z_3) &=&  
\left(
\begin{array}{cccc}
z_2-z_1 & \frac {z_2^2-z_1^2} 2 & z_3-z_1 & \frac {z_3^2-z_1^2} 2\\
 \frac {z_2^2-z_1^2} 2 &  \frac {z_2^3-z_1^3} 3  & \frac {z_3^2-z_1^2} 2&  \frac {z_3^3-z_1^3} 3\\
 \frac {z_2^3-z_1^3} 3 &  \frac {z_2^4-z_1^4} 4  & \frac {z_3^3-z_1^3} 3&  \frac {z_3^4-z_1^4} 4\\
 \frac {z_2^4-z_1^4} 4 &  \frac {z_2^5-z_1^5} 5  & \frac {z_3^4-z_1^4} 4&  \frac {z_3^5-z_1^5} 5
\end{array}
\right)
\\
\label{M2Diff}
\ \ M(0,\ z_2-z_1,\ z_3-z_1) &=& 
\left(
\begin{array}{cccc}
z_2-z_1 & \frac {(z_2-z_1)^2} 2 & z_3-z_1 & \frac {(z_3-z_1)^2} 2\\
 \frac {(z_2-z_1)^2} 2 &  \frac {(z_2-z_1)^3} 3  & \frac {(z_3-z_1)^2} 2&  \frac {(z_3-z_1)^3} 3\\
 \frac {(z_2-z_1)^3} 3 &  \frac {(z_2-z_1)^4} 4  & \frac {(z_3-z_1)^3} 3&  \frac {(z_3-z_1)^4} 4\\
 \frac {(z_2-z_1)^4} 4 &  \frac {(z_2-z_1)^5} 5  & \frac {(z_3-z_1)^4} 4&  \frac {(z_3-z_1)^5} 5
\end{array}
\right)
\\
\label{M2Norm}
M(0,z,1) &=& 
\left(
\begin{array}{cccc}
\vspace{0.3ex}
z & \frac {z^2} 2 & 1 & \frac 1 2\\
\vspace{0.3ex}
 \frac {z^2} 2 &  \frac {z^3} 3  & \frac 1 2&  \frac 1 3\\
\vspace{0.3ex}
 \frac {z^3} 3 &  \frac {z^4} 4  & \frac 1 3&  \frac 1 4\\
 \frac {z^4} 4 &  \frac {z^5} 5  & \frac 1 4&  \frac 1 5
\end{array}
\right)
\end{eqnarray*}

\section{A counter example in three space dimensions}

Let $T_0$ be the simplex of vertices
$
(0,0,0), \ (1,0,0), \ (0,1,0), \ (0,0,1)
$
and let $P_0$ be the harmonic polynomial of degree $4$
\begin{eqnarray*}
P_0 &:= & 3 x + 10 x^3 - 15 x^4 + 3 y - 18 x y - 15 x^2 y + 30 x^3 y - 15 x y^2 + 45 x^2 y^2 + 10 y^3 \\
& &+  30 x y^3 - 15 y^4 + 3 z - 18 x z - 15 x^2 z + 30 x^3 z - 18 y z + 240 x y z - 180 x^2 y z \\
& & - 15 y^2 z - 180 x y^2 z + 30 y^3 z - 15 x z^2 + 45 x^2 z^2 - 15 y z^2 - 180 x y z^2 + 45 y^2 z^2 \\
& & + 10 z^3 + 30 x z^3 + 30 y z^3 - 15 z^4.
\end{eqnarray*}
We define 
$$
u_0:=\nabla P_0 \in \cP_{2,3}.
$$
One can easily check using a formal computing program that all the linear functionals in $\cN_{2,3}(T_0)$ vanish on $u_0$, which shows that the conjecture Conj($2$,$3$), on quadratic vector fields in three dimensions, is not valid. 
The interested reader can download on the website arxiv.org, jointly to the preprint of this paper, a file that contains these verifications.

Finding a quadratic vector finite element for $\HCD$ in dimension $3$ thus remains an open question. Let us finally mention that, up to a multiplicative constant, $u_0$ is the only element of $\cP_{2,3}$ on which all the linear functionals $\cN_{2,3}(T_0)$ vanish.

\noindent\small
\nl
\nl
Jean-Marie Mirebeau
\nl
UPMC Univ Paris 06, UMR 7598, Laboratoire Jacques-Louis Lions, F-75005, Paris, France
\nl
CNRS, UMR 7598, Laboratoire Jacques-Louis Lions, F-75005, Paris, France
\nl
mirebeau@ann.jussieu.fr

\end{document}